\documentclass[12pt]{amsart}

\usepackage{color}
\definecolor{darkgreen}{rgb}{0, 0.40, 0}

\usepackage[
ps2pdf,                
colorlinks=true,
linkcolor=blue,
citecolor=darkgreen,
urlcolor=magenta
]{hyperref}

\usepackage{epsfig}
\usepackage{latexsym}
\usepackage{psfrag}
\usepackage{amssymb}
\usepackage{amsmath}
\usepackage{amsthm} 
\usepackage{verbatim} 

\newcommand{\calC}{{\mathcal{C}}}

\newcommand{\calF}{{\mathcal{F}}}

\newcommand{\NN}{{\mathbb{N}}}

\newcommand{\ZZ}{{\mathbb{Z}}}

\renewcommand{\setminus}{{\smallsetminus}}

\newcommand{\Id}{{\operatorname{Id}}}

\newcommand{\st}{{\:\mid\:}}
\newcommand{\from}{{\colon}} 

\def\defeq{\mathbin{:=}}  

\newcommand{\homeo}{{\medspace \cong \medspace}} 
 
\newcommand{\closure}[1]{{\overline{#1}}}

\newcommand{\neigh}{{\operatorname{neigh}}} 
 
\newcommand{\bdy}{{\partial}} 

\newcommand{\diam}{{\operatorname{diam}}} 

\newcommand{\MCG}{{\mathcal{MCG}}} 

\newcommand{\Teich}{{Teichm\"uller~}}

\newsavebox{\savepar}

\theoremstyle{plain}
\newtheorem{theorem}{Theorem}[section]

\newtheorem{lemma}[theorem]{Lemma}

\newtheorem{proposition}[theorem]{Proposition}

\theoremstyle{definition}

\newtheorem{remark}[theorem]{Remark}

\newtheorem{question}[theorem]{Question}

\theoremstyle{plain}
\newtheorem*{Thm:ComplementOfABall}{Theorem~\ref{Thm:ComplementOfABall}}

\newcommand{\dS}{{\dot{S}}}
\newcommand{\dX}{{\dot{X}}}
\newcommand{\dB}{{\dot{B}}}

\begin{document}

\title{The end of the curve complex}

\author{Saul Schleimer}
\address{\hskip-\parindent
        Department of Mathematics\\
	Rutgers University\\
        Piscataway, New Jersey 08854}
\email{saulsch@math.rutgers.edu}
\urladdr{ http://www.math.rutgers.edu/$\sim$saulsch}

\thanks{This work is in the public domain.}

\date{\today}

\begin{abstract}
Suppose that $S$ is a surface of genus two or more, with exactly one
boundary component.  Then the curve complex of $S$ has one end.
\end{abstract}

\maketitle

\section{Introduction}

We denote the compact, connected, orientable surface of genus $g$ with
$b$ boundary components by $S_{g,b}$.  The {\em complexity} of $S =
S_{g,b}$ is $\zeta(S) \defeq 3g - 3 + b$.  A simple closed curve $\alpha$
in $S$ is {\em essential} if $\alpha$ does not cut a disk out of $S$.
Also, $\alpha$ is {\em non-peripheral} if it does not cut an annulus
out of $S$.

When $\zeta(S) \geq 2$ the {\em complex of curves}, $\calC(S)$, is the
simplicial complex where vertices are isotopy classes of essential
non-peripheral curves.  The $k$--simplices are collections of $k+1$
distinct vertices having disjoint representatives.  We regard every
simplex as a Euclidean simplex of side-length one. If $\alpha$ and
$\beta$ are vertices of $\calC(S)$ let $d_S(\alpha, \beta)$ denote the
distance between $\alpha$ and $\beta$ in the one-skeleton
$\calC^1(S)$.  It is a pleasant exercise to prove that $\calC(S)$ is
connected.  It is an important theorem of H.\,Masur and
Y.\,Minsky~\cite{MasurMinsky99} that $\calC(S)$ is Gromov hyperbolic.

Let $B(\omega, r) \defeq \{ \alpha \in \calC^0(S) \st d_S(\alpha, \omega)
\leq r \}$ be the ball of radius $r$ about the vertex $\omega$.  We
will prove:

\begin{Thm:ComplementOfABall}
Fix $S \defeq S_{g,1}$ for some $g \geq 2$.  For any vertex $\omega
\in C(S)$ and for any $r \in \NN$: the complex spanned by $\calC^0(S)
\setminus B(\omega, r)$ is connected.
\end{Thm:ComplementOfABall}

For such surfaces, Theorem~\ref{Thm:ComplementOfABall} directly
answers a question of Masur's.  It also answers a question of G.\,Bell
and K.\,Fujiwara~\cite{BellFujiwara05} in the negative: the complex of
curves need not be quasi-isometric to a tree.
Theorem~\ref{Thm:ComplementOfABall} is also evidence for a positive
answer to a question of P.\,Storm:

\begin{question}
\label{Que:Storm}
Is the Gromov boundary of $\calC(S)$ connected?
\end{question}

Note that Theorem~\ref{Thm:ComplementOfABall} is only evidence for,
and not an answer to, Storm's question: for example, there is a
one-ended hyperbolic space where the Gromov boundary is a pair of
points.
Finally, as we shall see in Remark~\ref{Rem:NoQuotient}, it is not
obvious how to generalize Theorem~\ref{Thm:ComplementOfABall} to
surfaces with more (or fewer) boundary components.

\subsection*{Acknowledgments} 

I thank both Jason Behrstock and Christopher Leininger; their
observations, here recorded as Proposition~\ref{Prop:FibreIsBent},
were the origin of my thoughts leading to
Theorem~\ref{Thm:ComplementOfABall}.  I also thank Kenneth Bromberg
for showing me a simplification, given below, of my original proof of
Theorem~\ref{Thm:ComplementOfABall}.  Finally, I am grateful to Howard
Masur for both posing the main question and for many enlightening
conversations.

\section{Definitions and necessary results}

An important point elided above is how to define $\calC(S)$ when
$\zeta(S) = 1$.  The complex as defined is disconnected in these
cases.  Instead we allow a $k$--simplex to be a collection of $k+1$
distinct vertices which have representatives with small intersection.
For $S_{1,1}$ exactly one intersection point is allowed while
$S_{0,4}$ requires two.  In both cases $\calC(S)$ is the famous {\em
Farey tessellation}.  Note that $\calC(S_{0,3})$ is empty.  We will
not need to consider the other low complexity surfaces: the sphere,
the disk, the annulus, and the torus.

A subsurface $X \subset S$ is {\em essential} if every component of
$\bdy X$ is essential in $S$.  We will generally assume that $\zeta(X)
\geq 1$.  A pair of curves, or a curve and a subsurface, are {\em
tight} if they cannot be isotoped to reduce intersection.  We will
generally assume that all curves and subsurfaces discussed are tight
with respect to each other. We say a curve $\alpha$ {\em cuts} $X$ if
$\alpha \cap X \neq \emptyset$.  If $\alpha \cap X = \emptyset$ then
we say $\alpha$ {\em misses} $X$.

Following Masur and Minsky~\cite{MasurMinsky00}, we define the {\em
subsurface projection} map $\pi_X$: this maps vertices of $\calC(S)$
to collections of vertices of $\calC(X)$.  Fix a vertex $\alpha \in
\calC(S)$ and, for every component $\delta \subset \alpha \cap X$,
form $N_{\delta} \defeq \closure{\neigh(\delta \cup \bdy X)}$, a
closed regular neighborhood of $\delta \cup \bdy X$.  Take
$\pi_X(\alpha)$ to be the set of all vertices of $\calC(X)$ which
appear as a boundary component of some $N_{\delta}$.  If $\alpha$
misses $X$ then $\pi_X(\alpha) = \emptyset$.  Note if $\alpha \subset
S$ is contained in $X$ after tightening then $\pi_X(\alpha) = \{
\alpha \}$.

As a useful bit of notation, if $\alpha$ and $\beta$ both cut $X$, we
set
$$d_X(\alpha, \beta) \defeq \diam_X(\pi_X(\alpha), \pi_X(\beta))$$ with
diameter computed in $\calC^1(X)$.  Masur and Minsky give an
combinatorial proof~\cite[Lemma~2.2]{MasurMinsky00} that:

\begin{lemma}
\label{Lem:Disjoint}
If $\alpha$ and $\beta$ both cut $X$ and $d_S(\alpha, \beta) \leq 1$
then $d_X(\alpha, \beta) \leq 2$. \qed
\end{lemma}

By {\em geodesic} in $\calC(S)$ we will always be referring to a
geodesic in the one-skeleton.  Since $\calC(S)$ is Gromov hyperbolic
the exact position of the geodesic is irrelevant; we often use the
notation $[\alpha, \beta]$ as if the geodesic was determined by its
endpoints.  We immediately deduce from Lemma~\ref{Lem:Disjoint}:

\begin{lemma}
\label{Lem:WeakBoundedImage}
Suppose that $\alpha, \beta$ are vertices of $\calC(S)$, both cutting
$X$.  Suppose that $d_X(\alpha, \beta) > 2 \cdot d_S(\alpha, \beta)$.
Then every geodesic $[\alpha, \beta] \subset \calC(S)$ has a vertex
which misses $X$.  \qed
\end{lemma}

This is essentially Lemma~2.3 of~\cite{MasurMinsky00}.


\begin{remark}
\label{Rem:UniqueMissing}
There is a useful special case of Lemma~\ref{Lem:WeakBoundedImage}:
assume all the hypotheses and in addition that $\gamma$ is the unique
vertex of $\calC(S)$ missing $X$.  Then every geodesic connecting
$\alpha$ to $\beta$ contains $\gamma$.

In fact, $\gamma$ is the unique vertex missing $X$ exactly when $S
\setminus \neigh(\gamma) = X$ or $S \setminus \neigh(\gamma) = X \cup
P$ with $P \homeo S_{0,3}$: a {\em pants}.
\end{remark}

\begin{remark}
\label{Rem:NeedForBoundedImage}
Note that Lemma~\ref{Lem:WeakBoundedImage} is a weak form of the
Bounded Geodesic Image Theorem~\cite[Theorem~3.1]{MasurMinsky00}.  The
proof of their stronger result appears to require techniques from
\Teich theory.
\end{remark}

We now turn to the {\em mapping class group} $\MCG(S)$: the group of
isotopy classes of homeomorphisms of $S$.  Note that the natural
action of $\MCG(S)$ on $\calC(S)$ is via isometries.  We have an
important fact:

\begin{lemma}
\label{Lem:PA}
If $\psi \from S \to S$ is a pseudo-Anosov and $\alpha$ is a vertex of
$\calC^0(S)$ then $\diam_S(\psi^n(\alpha) \st n \in \ZZ)$ is
infinite. \qed
\end{lemma}

It follows that the diameter of $\calC(S)$ is infinite whenever
$\zeta(S) \geq 1$.  A proof of Lemma~\ref{Lem:PA}, relying on
Kobayashi's paper~\cite{Kobayashi88b}, may be found in the remarks
following Lemma~4.6 of~\cite{MasurMinsky00}.  As a matter of fact,
Masur and Minsky there prove more using train track machinery: any
orbit of a pseudo-Anosov map is a quasi-geodesic.  We will not need
this sharper version.

Note that if $\psi \from S \to S$ is a homeomorphism then we may
restrict $\psi$ to the curve complex of a subsurface $\psi|X \from
\calC(X) \to \calC(\psi(X))$.  This restriction behaves well with
respect to subsurface projection: that is, $\pi_{\psi(X)} \circ \psi =
\psi|X \circ \pi_X $.

We conclude this discussion by examining {\em partial maps}.  Suppose
that $X \subset S$ is an essential surface, not homeomorphic to $S$.
If $\psi \from S \to S$ has the property that $\psi|\closure{S
\setminus X} = \Id|\closure{S \setminus X}$ then we call $\psi$ a {\em
partial map} supported on $X$.  Note that if $\psi$ is supported on
$X$ then the orbits of $\psi$ do not have infinite diameter in
$\calC(S)$.  Since $\psi$ fixes $\bdy X$ and acts on $\calC(S)$ via
isometry, every point of an orbit has the same distance to $\bdy X$ in
$\calC(S)$.  Nonetheless, Lemmas~\ref{Lem:WeakBoundedImage}
and~\ref{Lem:PA} imply:

\begin{lemma}
\label{Lem:Swing}
Suppose $\psi \from S \to S$ is supported on $X$ and $\psi|X$ is
pseudo-Anosov.  Fix a vertex $\sigma \in \calC(S)$ and define $\sigma_n
\defeq \psi^n(\sigma)$.  Then for any $K \in \NN$ there is a power $n \in
\ZZ$ so that $d_X(\sigma, \sigma_n) \geq K$.  In particular, if $K > 4
\cdot d_S(\sigma, \bdy X)$ then every geodesic $[\sigma, \sigma_n]
\subset \calC(S)$ contains a vertex which misses $X$.  \qed
\end{lemma}

\section{No dead ends}

We require a pair of tools in order to prove
Theorem~\ref{Thm:ComplementOfABall}.  The first is:

\begin{proposition}
\label{Prop:NoDeadEnds}
Fix $S = S_{g,b}$.  For any vertex $\omega \in \calC(S)$ and for any
$r \in \NN$: every component of the subcomplex spanned by $\calC^0(S)
\setminus B(\omega, r)$ has infinite diameter.
\end{proposition}

A more pithy phrasing might be: the complex of curves has no {\em dead
ends}.  Proposition~\ref{Prop:NoDeadEnds} allows us to push vertices
away from $\omega$ while remaining inside the same component of
$\calC(S) \setminus B(\omega, r)$.  The proof is a bit subtle due to
the behavior of $\calC(S)$ near a non-separating curve.

\begin{proof}[Proof of Proposition~\ref{Prop:NoDeadEnds}]
If $S = S_{0,3}$ is a pants then the curve complex is empty and there
is noting to prove.  If $\calC(S)$ is a copy of the Farey graph then
the claim is an easy exercise.  
So we may suppose that $\zeta(S) \geq 2$.

Now fix a vertex $\alpha \in \calC(S) \setminus B(\omega, r)$.  Set $n
\defeq d_S(\alpha, \omega)$.  Thus $n > r$.  Our goal is to find a curve
$\delta$, connected to $\alpha$ in the complement of $B(\omega, n -
1)$, with $d_S(\delta, \omega) = n + 1$.  Doing this repeatedly proves
the proposition.  Note that finding such a vertex $\delta$ is
straight-forward if $r = 0$ and $n = 1$.  This is because $\calC(S)
\setminus \omega$ is connected and because, following
Lemma~\ref{Lem:PA}, we know that the diameter of $\calC(S)$ is
infinite.  Henceforth we will assume that $n \geq 2$; that is,
$\omega$ cuts $\alpha$.

Fix attention on a component $X$ of $S \setminus \neigh(\alpha)$ which
is not a pair of pants.  So $\zeta(X) \geq 1$ and, by the comments
following Lemma~\ref{Lem:PA}, $\calC(X)$ has infinite diameter.  Since
$\omega$ cuts $\alpha$ we find that $\omega$ also cuts $X$.  Choose a
curve $\beta$ contained in $X$ with $d_X(\beta, \omega) \geq 2n + 1$.
Note that $d_S(\alpha, \beta) = 1$.  We may assume that $\beta$ is
either non-separating or cuts a pants off of $S$.  (To see this: if
$\beta$ cannot be chosen to be non-separating then $X$ is planar.  As
$\zeta(X) \geq 1$ we deduce that $X$ has at least four boundary
components.  At most two of these are parallel to $\alpha$.)  It
follows from Lemma~\ref{Lem:WeakBoundedImage} that any geodesic from
$\beta$ to $\omega$ in $\calC(S)$ has a vertex $\gamma$ which misses
$X$.

By the triangle inequality $d_S(\gamma, \omega)$ equals $n$ or $n -
1$.  In the former case we are done: simply take $\delta = \beta$ and
notice that $d_S(\beta, \omega) = n + 1$.  In the latter case
$d_S(\beta, \omega) = n$ and we proceed as follows: replace $\alpha$
by $\beta$ and replace $X$ by $Z \defeq S \setminus \neigh(\beta)$.
We may now choose $\delta$ to be a vertex of $\calC(Z)$ with
$d_Z(\delta, \omega) \geq 2n + 1$.  As above, any geodesic $[\delta,
\omega] \subset \calC(S)$ has a vertex which misses $Z$.  Since
$\beta$ is the {\em unique} vertex not cutting $Z$
Remark~\ref{Rem:UniqueMissing} implies that $\beta \in [\delta,
\omega]$.  Thus $d_S(\delta, \omega) = n + 1$ and we are done.
\end{proof}

\section{The Birman short exact sequence}

We now discuss the second tool needed in the proof of
Theorem~\ref{Thm:ComplementOfABall}.  Following Kra's notation
in~\cite{Kra81} let $\dS = S_{g,1}$ and $S = S_g$ for a fixed $g \geq
2$.  Let $\rho \from \dS \to S$ be the quotient map crushing $\bdy
\dS$ to a point, say $x \in S$.  This leads to the {\em Birman short
exact sequence}:
$$\pi_1(S, x) \to \MCG(\dS) \to \MCG(S)$$ for $g \geq 2$.  The map
$\rho$ gives the second arrow.  The first arrow is defined by sending
$\gamma \in \pi_1(S, x_0)$ to a mapping class $\psi_\gamma$.  There is
a representative of this class which is isotopic to the identity, in
$S$, via an isotopy dragging $x$ along the path $\gamma$.  See
Birman's book~\cite{Birman74} or Kra's paper~\cite{Kra81} for further
details.

Fix an essential subsurface $\dX \subset \dS$ and let $X =
\rho(\dX)$.  If $\gamma \in \pi_1(S, x)$ is contained in $X$ then
$\psi_\gamma$ is a partial map, supported in $\dX$.  We say that
$\gamma$ {\em fills} $X$ if $\gamma \subset X$ and, in addition, every
representative of the free homotopy class of $\gamma$ cuts $X$ into a
collection of disks and peripheral annuli.  For future use we record a
well-known theorem of I.\,Kra~\cite{Kra81}:

\begin{theorem}
\label{Thm:Kra}
Suppose that $\zeta(\dX) \geq 1$.  If $\gamma$ fills $X$ then
$\psi_\gamma|\dX$ is pseudo-Anosov. \qed
\end{theorem}


Now note that, corresponding to the Birman short exact sequence, there
is a ``fibre bundle'' of curve complexes:
$$\calF_\tau \to \calC(\dS) \to \calC(S).$$ Here $\tau$ is an
arbitrary vertex of $\calC(S)$ and $\calF_\tau \defeq
\rho^{-1}(\tau)$.  The second arrow is given by $\rho$.  The first is
the inclusion of $\calF_\tau$ into $\calC(\dS)$.


\begin{remark}
\label{Rem:NoQuotient}
If $|\bdy S| \geq 2$ then collapsing one boundary component does not
induce a map on the associated curve complexes.  Thus, it is not clear
how to generalize Theorem~\ref{Thm:ComplementOfABall} to such
surfaces.  If $\bdy S$ is empty then I do not know of any
interesting quotients or electrifications of $\calC(S)$.
\end{remark}

Using the Birman short exact sequence we obtain an action of $\pi_1(S,
x)$ on the curve complex $\calC(\dS)$.  Behrstock and Leininger observe
that:

\begin{proposition}
\label{Prop:FibreIsBent}
The map $\rho \from \calC(\dS) \to \calC(S)$ has the following properties:
\begin{itemize}
\item
It is $1$--Lipschitz.
\item
For any $\alpha \in \calC(\dS)$, $\gamma \in \pi_1(S, x)$ we have
$\rho(\alpha) = \rho(\psi_\gamma(\alpha))$.
\item
Every fibre $\calF_\tau$ is connected.
\end{itemize}
\end{proposition}

\begin{remark}
\label{Rem:FibreIsBent}
Behrstock and Leininger's interest in the fibre $\calF_\tau$ was to
give a ``natural'' subcomplex of $\calC(S)$ which is not quasi-convex:
this is implied by the first pair of properties.
\end{remark}

\begin{remark}
\label{Rem:FibreIsBentTree}
More of the structure of $\calF_\tau$ is known. For example, since $S$
is closed, the fibre $\calF_\tau$ is either a single $\pi_1(S,
x)$--orbit or the union of a pair of orbits depending on whether
$\tau$ is non-separating or separating.  Furthermore, $\calF_\tau$ is
a tree.  See~\cite{KentEtAl06} for a detailed discussion.
\end{remark}

\begin{proof}[Proof of Proposition~\ref{Prop:FibreIsBent}]
Fix an essential non-peripheral curve $\alpha$ in $\dS$.  Note that
$\rho(\alpha)$ is essential in $S$ and so the induced map $\rho \from
\calC(\dS) \to \calC(S)$ is well-defined.  If $\alpha$ and $\beta$ are
disjoint in $\dS$ then so are their images in $S$.  Thus $\rho$ does
not increase distance between vertices and the first conclusion holds.


Now fix a curve $\alpha \subset \dS$ and $\gamma \in \pi_1(S,x)$.
Note that $\psi_\gamma$ is isotopic to the identity in $S$.  Thus the
images $\rho(\psi_\gamma(\alpha))$ and $\rho(\alpha)$ are isotopic
in $S$.  It follows that $\rho(\alpha) = \rho(\psi_\gamma(\alpha))$
as vertices of $\calC(S)$, as desired.

Finally, fix $\tau \in \calC(S)$.  Let $\calF_\tau$ be the fibre over
$\tau$.  Pick $\alpha, \beta \in \calF_\tau$.  It follows that $a
\defeq \rho(\alpha)$ and $b \defeq \rho(\beta)$ are both isotopic to
$\tau$ and so to each other.  We induct on the intersection number
$\iota(\alpha, \beta)$.  Suppose the intersection number is zero.
Then $\alpha$ and $\beta$ are disjoint and we are done.  Suppose that
the intersection number is non-zero.  Since $a$ and $b$ are isotopic,
yet intersect, they are {\em not} tight with respect to each other.
It follows that there is a bigon $B \subset S \setminus (a \cup b)$.
Since $\alpha$ and $\beta$ are tight in $\dS$ the point $x$ must lie
in $B$.  Let $\dB \defeq \rho^{-1}(\bar{B})$.  Now construct a curve
$\beta' \subset \dS$ by starting with $\beta$, deleting the arc $\beta
\cap \dB$, and adding the arc $\alpha \cap \dB$.  Isotope $\beta'$ to
be tight with respect to $\alpha$.  Now $\beta' \in \calF_\tau$
because $\rho(\beta')$ is isotopic to $\rho(\beta)$ in $S$.  Finally,
$\iota(\alpha, \beta') \leq \iota(\alpha, \beta) - 2$.
\end{proof}

\section{Proving the theorem}

We are now equipped to prove: 

\begin{theorem}
\label{Thm:ComplementOfABall}
Fix $\dS \defeq S_{g,1}$ for some $g \geq 2$.  For any vertex $\omega \in
C(\dS)$ and for any $r \in \NN$: the complex spanned by $\calC^0(\dS)
\setminus B(\omega, r)$ is connected.
\end{theorem}

As above we use the notation $\dS = S_{g,1}$ and $S = S_g$ for some
fixed $g \geq 2$.  Also, we have defined a map $\rho \from \calC(\dS)
\to \calC(S)$ induced by collapsing $\bdy \dS$ to a point, $x$.  As
above we use $\calF_\tau = \rho^{-1}(\tau)$ to denote the fibre over
$\tau$.  

\begin{proof}[Proof of Theorem~\ref{Thm:ComplementOfABall}]
Choose $\alpha'$ and $\beta'$ vertices of $\calC(\dS) \setminus B(\omega,
r)$.  By Proposition~\ref{Prop:NoDeadEnds} we may connect $\alpha'$ and
$\beta'$, by paths disjoint from $B(\omega, r)$, to vertices outside of
$B(\omega, 3r)$.  Call these new vertices $\alpha$ and $\beta$.  We
may assume that both $\alpha$ and $\beta$ are non-separating because
such vertices are $1$--dense in $\calC(\dS)$.

Choose any vertex $\tau \in \calC(S)$ so that $d_S(\tau, \rho(\omega))
\geq 4r$.  This is always possible because $\calC(S)$ has infinite
diameter.  (See the remarks after Lemma~\ref{Lem:PA}.)  It follows
from Proposition~\ref{Prop:FibreIsBent} that $\calF_\tau \cap
B(\omega, r) = \emptyset$.  We will now connect each of $\alpha$ and
$\beta$ to some point of $\calF_\tau$ via a geodesic disjoint from
$B(\omega, r)$.  Since $\calF_\tau$ is connected, by
Proposition~\ref{Prop:FibreIsBent}, this will complete the proof of
Theorem~\ref{Thm:ComplementOfABall}.

Let $\dX \defeq \dS \setminus \alpha$ and take $X \defeq \rho(\dX)$.
Fix any point $\sigma$ in $\calF_\tau$.  If $\sigma = \alpha$ then
$\alpha$ is trivially connected to the fibre.  So suppose that
$\sigma \neq \alpha$.  Since $\alpha$ is non-separating deduce that
$\sigma$ cuts $\dX$.  Now, since $\zeta(\dS) \geq 4$ we have
$\zeta(\dX) \geq 3$.  Let $\gamma \in \pi_1(S, x)$ be any homotopy
class so that $\psi_\gamma$ is supported in $\dX$ and so that $\gamma$
fills $X$.  By Kra's Theorem (\ref{Thm:Kra}) $\psi_\gamma|\dX$ is
pseudo-Anosov.

Since $\calF_\tau$ is left setwise invariant by $\pi_1(S, x)$
(Proposition~\ref{Prop:FibreIsBent}) the curves $\sigma_n \defeq
\psi_\gamma^n(\sigma)$ all lie in $\calF_\tau$.  Since
$\psi_\gamma|\dX$ is pseudo-Anosov, Lemma~\ref{Lem:Swing} gives an $n
\in \ZZ$ so that every geodesic $g = [\sigma, \sigma_n] \subset
\calC(\dS)$ has a vertex which misses $\dX$.  Since $\alpha$ is
non-separating, as in Remark~\ref{Rem:UniqueMissing}, it follows that
$\alpha$ is actually a vertex of $g$.

We now claim that at least one of the two segments $[\sigma, \alpha]
\subset g$ or $[\alpha, \sigma_n] \subset g$ avoids the ball
$B(\omega, r)$.  For suppose not: then there are vertices $\mu, \mu'
\in g$ on opposite sides of $\alpha$ which both lie in $B(\omega,
r)$.  Thus $d_{\dS}(\mu, \mu') \leq 2r$.  Since $g$ is a geodesic the
length along $g$ between $\mu$ and $\mu'$ is at most $2r$.  Thus
$d_{\dS}(\omega, \alpha) \leq 2r$.  This is a contradiction.

Thus we can connect $\alpha$ to a vertex of $\calF_\tau$ (namely,
$\sigma$ or $\sigma_n$) avoiding $B(\omega, r)$.  Identically, we can
connect $\beta$ to a vertex of $\calF_\tau$ while avoiding $B(\omega,
r)$.  As noted above, this completes the proof.
\end{proof}

\bibliographystyle{plain}
\bibliography{bibfile}
\end{document}